\documentclass[conference,a4paper]{IEEEtran}
\IEEEoverridecommandlockouts
\usepackage{cite}
\usepackage{amsmath,amssymb,amsfonts}
\usepackage{algorithmic}
\usepackage{graphicx}
\usepackage{subcaption}
\usepackage{textcomp}
\usepackage{xcolor}
\usepackage{stackengine}
\usepackage{color}

\def\BibTeX{{\rm B\kern-.05em{\sc i\kern-.025em b}\kern-.08em
  T\kern-.1667em\lower.7ex\hbox{E}\kern-.125emX}}
 
\begin{document}

\title{The Energy-Delay Pareto Front in Cache-enabled Integrated Access and Backhaul mmWave HetNets\\}
\author{\parbox{6 in}{\centering Wen Shang and Vasilis Friderikos \\
         Department of Engineering\\
        King's College London, London WC2R 2LS, U.K.\\
        E-mail:\{wen.shang, vasilis.friderikos\}@kcl.ac.uk}
}

\maketitle

\begin{abstract}
In this paper, to address backhaul capacity bottleneck and concurrently optimize energy consumption and delay, we formulate a novel weighted-sum multi-objective optimization problem where popular content caching placement and integrated access and backhaul (IAB) millimeter (mmWave) bandwidth partitioning are optimized jointly to provide Pareto efficient optimal non-dominating solutions. In such integrated networks analysis of what-if scenarios to understand trade-offs in decision space, without losing sight of optimality, is important. A wide set of numerical investigations reveal that compared with the nominal single objective optimization schemes such as optimizing only the delay or the energy consumption the proposed optimization framework allows for a reduction of the aggregation of energy consumption and delay by an average of $30\%$ to $55\%$.

\end{abstract}
\begin{IEEEkeywords}
Beyond 5G, mmWave communications, wireless communications, network optimization, Pareto optimization,  Integrated Access and Backhaul.
\end{IEEEkeywords}

\section{Introduction}
\IEEEPARstart{T}he deployment of small cell base stations (SBSs) overlaid with macro cell base station (MBS) to construct a heterogeneous network (HetNet) has been identified as one of the most effective technologies for supporting the envisioned explosive traffic loads in fifth generation (5G) and beyond (B5G) cellular networks\cite{lv2017interference}. A key challenge, however, is that as the network become denser, the limited capacity of the backhaul links is becoming a bottleneck, and as such limiting the capabilities of the network. Although equipping each SBS with a wired backhaul link can alleviate the aforementioned problem and provide more strict guarantees on the performance \cite{bacstuug2016delay}, it comes at a high infrastructure cost. On the other hand, connecting the SBSs to the core network via traditional microwave is a low-cost option, but it may not meet throughput requirements especially in densely deployed small cells \cite{chandrasekhar2008femtocell}.

In fact, it has been demonstrated that, notable data traffic is due to multiple requests for the same content\cite{khan2019energy}. Popular files can be cached at small cells to offload data traffic to the edge(\cite{bacstuug2016delay}, \cite{khan2019energy, wu2020energy, song2021energy, gabry2016energy}) resulting in alleviating traffic load on backhaul links while meeting link capacity and throughput requirements. In a closely related work, Wu and Lu in \cite{wu2020energy} investigated the trade-off between energy consumption and delay in cache-enabled small cell networks, analysing power allocation and content placement in the edge to increase network efficiency, but only on single-tier networks. Two-tier network has been considered in \cite{song2021energy}. Energy harvesting and collaboration between SBSs have been proposed to reduce energy consumption and traffic load on the backhaul links. Because of the limited cache capacity, network densification is still restricted by the traditional microwave backhaul link. The work in \cite{gabry2016energy} and\cite{lv2017interference} investigate energy consumption and throughput optimization, respectively, but focusing on one target may not improve overall network efficiency.

Apart from edge caching, another effective way to tackle backhaul link capacity bottleneck is to exploit new spectrum bands, such as millimeter wave (mmWave) spectrum. The work in \cite{wang2019energy} assessed energy efficiency and delay in mmWave networks, but edge caching was not taken explicitly into account. Capitalizing on the availability of the bandwidth spectrum of mmWave, it is possible to establish access and backhaul links by using the same bandwidth and thus share the same infrastructure to reduce cost expenditure without penalizing performance. This is the core idea behind the so-called integrated access and backhaul (IAB) concept that has been proposed by the Third Generation Partnership Project (3GPP) \cite{zhang2021survey}. 
The effect of bandwidth partitioning between access and backhaul links on the achievable data rate has been observed in \cite{saha2019millimeter}, but it has been pointed out that offloading users to edge may not provide similar improvements in an IAB system as it would be in a fiber-backhauled network. This is an intuitively expected result owing to the fact that without implementing efficient edge-caching, data traffic offloaded to SBSs actually returns to the MBS via the self-backhaul links.

While the aforementioned works have clearly demonstrated the benefits that an mmWave architecture can bring, the trade-off between energy consumption and delay has not been meticulously analyzed in a cache-enabled HetNet where IAB mmWave is used to address backhaul bottleneck. To fill the above void in this work we explicitly consider the trade-off via a multi-objective optimization problem where cache placement and bandwidth allocation are jointly considered under a IAB mmWave spectrum framework. To this end, the main contributions of the paper are described as follows:

\begin{itemize}
  \item a novel multi-objective optimization framework is proposed under cache-enabled IAB mmWave HetNets to reveal the Pareto efficient frontier of the underlying energy-delay trade-off;
  \item a dynamic bandwidth partition and edge caching allocation are considered jointly that allows a decision maker to create different operation modes- all of them being Pareto optimal; 
  \item we also shed light on the impact on the performance of the caching time, an aspect which hasn't been previously reported, in conjunction with the small cell density, the available cache capacity and transmission power;
  \item results of delay-only, energy-only and no-edge-caching schemes in according to the structures in \cite{bacstuug2016delay}, \cite{gabry2016energy}, \cite{saha2019millimeter}, respectively, are investigated as the baseline and compared with the proposed framework.
\end{itemize}

The rest of the paper is organized as follows. The system model is described in Section \ref{introduction}. Then, in Section \ref{problem}, the multi-objective minimization problem is formulated. Numerical investigations over key metrics are presented in Section \ref{results}. Finally, concluding remarks are discussed in Section \ref{conclusion}.

\section{System Model}\label{introduction}
In this paper, a two-tier cache-enabled HetNet is considered, similar to \cite{song2021energy} with the key difference that a IAB framework based on mmWave is used. To obtain optimal popular content caching decision and bandwidth allocation, energy consumption and delay are optimized jointly via a proposed weighted-sum optimization problem.
\subsection{IAB mmWave Heterogeneous Network Topology}

\begin{figure}[htbp]
\centerline{\includegraphics[width=0.5\textwidth]{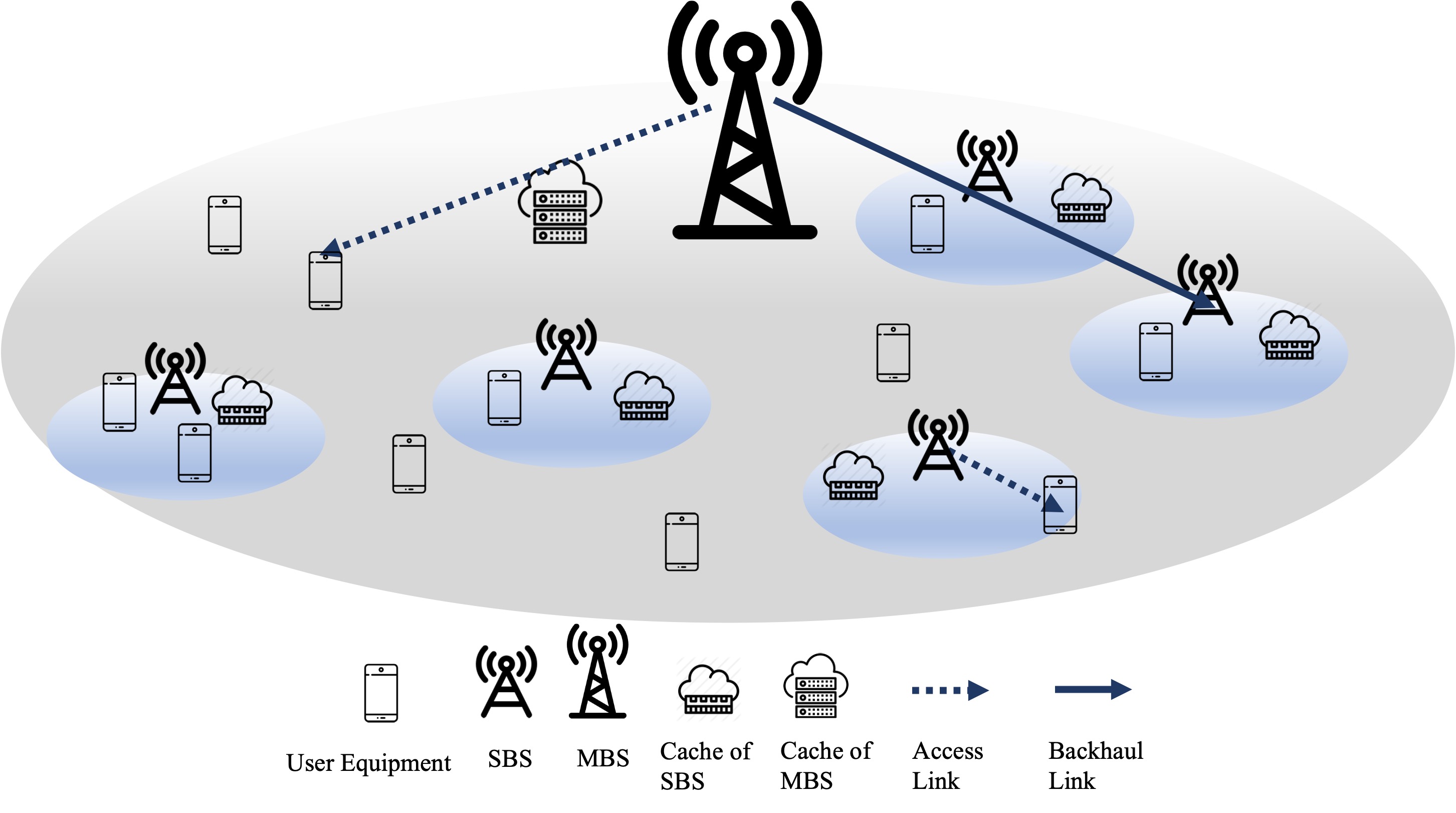}}
\caption{The IAB millimeter (mmWave) heterogeneous networks network topology considered in this work.}
\label{fig:Network Topology}
\end{figure}

As shown in Fig. \ref{fig:Network Topology}, the network topology under investigation comprises one MBS which is fiber backhauled to a core network data center where popular contents requested by users can be stored, $B$ SBSs with limited cache storage and $U$ user equipments (UEs). We denote all BSs by $\mathcal{B}=\{0,1,\cdots,B\}$, where $j=0$ represents the MBS located in the center and $j=1, \cdots,B$ represents the $j^{th}$ SBS randomly distributed around the MBS. The set of UEs, $\mathcal{U}=\{1,\cdots,U\}$, are also randomly distributed and associated with one BS. We assume that during off-peak time, the MBS offloads popular content to the cache memory of SBSs .

\subsection{Network Caching Model}
In the aforementioned network system, a popular content set requested by UEs, is cached at the core data center and indexed by popularity. Users request files according to an a-priori known file popularity distribution, and we assume the popularity of files, $p_m$, $m\in \mathcal{K}=\{1,2\cdots,K\}$, follows a widely used Zipf-like distribution \cite{gabry2016energy},
\[p_{m}=\frac{m^{-r}}{\sum_{k\in\mathcal{K}}k^{-r}}\]
where $r$ is the distribution's skewness coefficient. Furthermore, size of each file is assumed to be $Q$ bits, and the cache capacity of each SBS equals to $N\cdot Q$ bits, which indicating that it can cache $N$ files at most.

The MBS has a coverage area of radius $R_M$, while each SBS has a smaller coverage of radius $R_S$ overlaid in the macro cell. With $a_{ij}(i\in\mathcal{U},j\in\mathcal{B}),$ we refer to the association relationship between UEs and BSs, and each UE is associated with only one BS. The $i^{th}$ UE located within the coverage of SBSs, will be associated with the nearest $j^{th}$ SBS and $a_{ij}=1$. Files requested by the user which have been offloaded to the edge, will be retrieved from the SBS directly via access link; in that case we have a cache hit. Otherwise, the request will be redirected to the MBS, this case is denoted as a cache miss. We note that UEs located out of coverage of all SBSs will be connected to the MBS directly, i.e., $a_{i0}=1$.

\subsection{Spectrum Allocation Model for mmWave IAB}

As stated in \cite{lv2017interference}, access and backhaul links can employ distinct or same portion of the available spectrum, i.e., in-band or out-of-band mode. Similar to \cite{saha2019millimeter}, the out-of-band mode is applied to analyze the IAB mmWave HetNet performance. This is because we can capitalize on the significant potential of supporting large bandwidth in the mmWave bands and allow access and backhaul links to operate in different portions of that spectrum to limit cross-tire interference. As shown in \cite{saha2019millimeter}, a proper spectrum resource split between access and backhaul links can have a significant impact on the system performance.

\begin{figure}[htbp]
\centering
\includegraphics[width=0.3\textwidth]{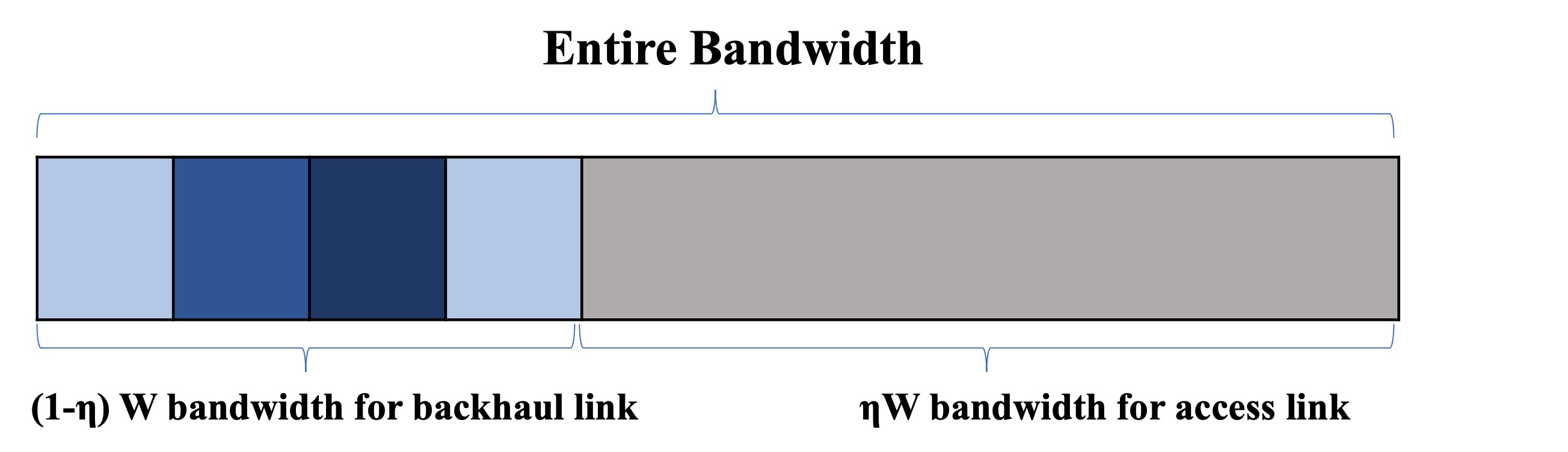}
\caption{\label{fig:bandwidth allocation} Graphical representation of the spectrum allocation model for mmWave IAB.}
\end{figure}

\subsection{Wireless Communication Model}
As shown in Fig. \ref{fig:bandwidth allocation}, we assume that a fraction of $\eta  (0\le \eta \le 1) $ bandwidth is allocated to the access links shared by all BSs with frequency reuse and the remaining is allocated in an orthogonal manner to backhaul links, i.e., wireless links between the MBS and the SBSs. 
When the $i^{th}$ UE is connected to the $j^{th}$ BS retrieving the $m^{th}$ content, the Signal to Interference plus Noise ratio (SINR), $SINR_{ij}$, is given as,
\begin{equation}
    SINR_{ij}=\frac{P_{j}^t|h_{ij}|^2G^a\beta ||d_{ij}||^{-\delta}}{N_0W+\sum_{\substack{j' \in \mathcal{J} \backslash \{j\}}}{I_{ij'}}}
\end{equation}

where $P_{j}^t$ denotes the SBS's transmission power, $|h_{ij}|^2$ refers to the small scale fading of the wireless communication link. Since Rayleigh fading assumption has been proven to be inaccurate for mmWave link communications, we assume an independent Nakagami fading for the transmission link \cite{bai2014coverage}, then $|h_{ij}|^2$ is a normalized Gamma random variable. Furthermore, $\beta=\frac{c}{(4\pi f_c)^2}$ is the propagation loss at a reference distance (1m) where $f_c$ is the carrier frequency. Also, $||d_{ij}||^{-\delta}$ expresses the path loss, where $||d_{ij}||$ denotes the distance between them. Note that $\delta$ is equal to $\delta^l$ for Line Of Sight (LOS) area and $\delta^n$ for Non-Line Of Sight (NLOS) area. Propagation loss in mmWave frequencies can be compensated by beam-formed directional transmission. Therefore, all BSs are assumed to be equipped with antenna arrays and $G^a$ refers to the antenna gain\cite{saha2019millimeter}. $I_{ij'}$ denotes the interference from the $j'^{th}$ BSs to the $i^{th}$ UE, and to provide bounds on the performance we consider the worst scenario where the UE receives interference from all other BSs. ${I_{ij'}}=P_{j'}G^a|h_{ij'}|^2\beta ||d_{ij'}||^{-\delta}$. Finally, $N_0W$ represents the noise power over the entire bandwidth $W$.

In cache miss, the request is redirected to the MBS, $SNR_{j0}$ refers to the Signal to Noise ratio (SNR) of the backhaul link between the MBS and the $j^{th}$ SBS, and can be described as follows,
\begin{equation}
    SNR_{j0}=\frac{P_{0}^t|h_{j0}|^2G^a\beta ||d_{j0}||^{-\delta}}{N_0W}
\end{equation}
in the above expression $P_0^t$ denotes the MBS transmission power, $|h_{j0}|^2$ and $||d_{j0}||^{-\delta}$ denotes the Nakagami fading and path loss between the MBS and the $j^{th}$ SBS, respectively.

Therefore, the achievable data rate of the access link between the $j^{th}$ BS and the $i^{th}$ UE can be written as follows,

\begin{equation}
 r_{ij}^a=\frac{\eta W}{\sum_{i\in\mathcal{U}} a_{ij}\sum_{m\in\mathcal{K}}p_{mi}}log_2(1+SINR_{ij}),\quad \forall j \in\mathcal{J}
 \label{con:3}
\end{equation}
where $\eta W$ is the bandwidth allocated to access link shared by all BSs, $p_{mi}$ refers to the probability of $i^{th}$ UE requesting the $m^{th}$ file, as stated before, $p_{mi}=p_m$,
and $\sum_{i\in\mathcal{U}} a_{ij}\sum_{m\in\mathcal{K}}p_{mi}$ is the corresponding traffic load, i.e. the aggregated number of requests. 

The achievable data rate of the backhaul link between MBS and $j^{th}$ SBS is given as follows,

\begin{equation}
\begin{aligned}
   r_{j0}^b=\frac{(1-\eta)W}{J\sum_{m\in\mathcal{K}}(1-x_{jm})\sum_{i\in\mathcal{U}}(p_{mi}a_{ij})} log_2(1+SNR_{j0}),&\\
   \quad \forall j \in\mathcal{J}\,and\,j\ne 0 &
\end{aligned}
\end{equation}

where $\frac{(1-\eta)W}{J}$ is the bandwidth allocated to backhaul link between MBS and $j^{th}$ SBS. The component $\sum_{m\in\mathcal{K}}(1-x_{jm})\sum_{i\in\mathcal{U}}(p_{mi}a_{ij})$ expresses the corresponding traffic load.

\subsection{Energy consumption and Delay}
\subsubsection{Energy Consumption}
We define the total energy consumption as
\begin{equation}
    E=CE+TE
\end{equation}
where $CE$ and $TE$ refer to content caching and wireless transmission energy consumption, respectively. Those two components are detailed in the sequel.

On one hand, the energy consumed by caching contents at the SBSs is given by\cite{xu2014coordinated},

\begin{equation}
\begin{aligned}
  CE&=\sum_{j\in\mathcal{B}\backslash 0} CE_j=\sum_{j\in\mathcal{B}\backslash 0}{\sum_{m\in\mathcal{K}}w^{ca}x_{jm}TQ}
\end{aligned}
\end{equation}
where $w^{ca}$ refers to the power efficiency of caching, $T$ is the caching time. With $x_{jm}$ we denote the caching decision of the $j^{th}$ SBS, where $x_{jm}=1$ denotes that the $m^{th}$ content is cached at the $j^{th}$ SBS, and $x_{jm}=0$ expresses a cache miss. 

On the other hand, energy consumed when delivering popular content to receivers consists of two parts, access link transmission consumption from BSs to UEs and additional backhaul link delivering when a cache miss occurs. Therefore, the wireless transmission energy consumption is defined as follows,

\begin{equation}
    TE=\sum_{j\in\mathcal{B}} TE_j^a+TE_0^b
\end{equation}
where $TE_j^a$ and $TE_0^b$ refer to the transmission consumption of the access link and the backhaul link, respectively. Similar to \cite{wu2020energy}, the mathematical characterization of the energy consumption can be written as follows,

\begin{equation}
    TE_j^a=\sum_{i\in\mathcal{U}}{P_{ij}^t({\frac{Q}{r_{ij}^a}\sum_{m\in\mathcal{K}}a_{ij}p_{mi}})}
\end{equation}
\begin{equation}
    TE_0^b=\sum_{j\in\mathcal{B}\backslash0}{P_{j0}^{t}{\frac{Q}{r_{j0}^b}}(\sum_{i\in\mathcal{U}}\sum_{m\in\mathcal{K}}a_{ij}p_{mi}}(1-x_{jm}))\label{con:9}
\end{equation}
where $P^{t}_{ij}$ and $P_{j0}^{t}$ refer to the transmitter transmission power allocated to a specific receiver, and fractional power control is used to efficiently distribute the transmission power to reduce energy consumption \cite{jindal2008fractional}. 

Therefore, from  equations (\ref{con:3})-(\ref{con:9}), the energy consumption can be written as follows,

\begin{equation}
    \begin{aligned}
    E&=CE+TE\\
    &=\sum_{j\in\mathcal{B}\backslash 0}CE_j+\sum_{j\in\mathcal{B}} TE_j^a+TE_0^b\\
    &=\sum_{j\in\mathcal{B}\backslash 0}\sum_{m\in\mathcal{K}}w^{ca}x_{jm}TQ+\\
    &\sum_{j\in\mathcal{B}}{\sum_{i\in\mathcal{U}}{P_{ij}^t( \frac{Q\sum_{i\in\mathcal{U}}{a_{ij}}}{\eta Wlog_2(1+SINR_{ij})}a_{ij})}}+\\
    &\sum_{j\in\mathcal{B}\backslash 0}{P_{j0}^t \frac{QJ(\sum_{m\in\mathcal{K}}{(\sum_{i\in\mathcal{U}}{p_{mi}a_{ij}})(1-x_{jm})})^2}{(1-\eta)Wlog_2(1+SNR_{j0})}}
    \end{aligned}
\end{equation}

\subsubsection{Aggregated Delay}

The delivery delay equals to the access wireless transmission time \cite{wu2020energy}, and an additional backhaul delay if cache miss occurs. Let $d_{ij}^m$ denotes the delay when the $i^{th}$ UE associated with the $j^{th}$ BS requests the $m^{th}$ file:

\begin{equation}
d_{ij}^m=\frac{Q}{r_{ij}^a}+\frac{Q}{r_{j0}^b}{(1-x_{jm})}
\end{equation}

 $\frac{Q}{r_{ij}^a}$ and $\frac{Q}{r_{j0}^b}$ refers to the access and backhaul link transmission delay. Therefore, the delay can be written as follows,
\begin{equation}
\begin{aligned}
    D&=\sum_{j\in\mathcal{B}}\sum_{i\in\mathcal{U}}\sum_{m\in\mathcal{K}}p_{mi}a_{ij}\cdot d_{ij}^m\\
    &=\sum_{j\in\mathcal{B}}\sum_{i\in\mathcal{U}}\sum_{m\in\mathcal{K}}p_{mi}a_{ij}(\frac{Q}{r_{ij}^a}+(1-x_{jm})\frac{Q}{r_{j0}^b})
    \end{aligned}
\end{equation}

\section{ Multi-Objective Optimization Problem}\label{problem}
As already eluded before, energy consumption and delay are both critical metrics in cache-enabled HetNets. In this section, we intend to optimize them jointly and reveal the trade-off between them via the Pareto efficient frontier. The inherent trade-off is captured by a weighted-sum multi-objective optimization problem \cite{marler2004survey} and we formulate it as follows,

\begin{equation}
    \min_{\substack{\mathbf{\eta,\mathbf{x}}}} \quad \alpha\delta_e E+(1-\alpha)\delta_d D
    \label{con:objective}
\end{equation}
\begin{subequations}
\begin{equation}
     s.t. \quad \,x_{jm}\in\{0,1\},\,\forall j\in\mathcal{B},\,\forall m\in\mathcal{K}\label{con:a}
\end{equation}
\begin{equation}
    \begin{aligned}
      \sum_{m\in\mathcal{K}} x_{jm}\le N, \,\forall j\in\mathcal{B}, j\ne 0 \label{con:b}
    \end{aligned}
\end{equation}
\begin{equation}
    \begin{aligned}
       0\le \eta \le 1 \label{con:c}
    \end{aligned}
\end{equation}
\begin{equation}
    r_{ij}a_{ij} \ge \gamma_i\label{con:d}
\end{equation}
\begin{equation}
    r_{j0}^b \ge \tau_j\label{con:e}
\end{equation}
\end{subequations}

In objective function, the weighting parameter $\alpha \in (0,1)$ can be altered to reflect the relative preferences over two objectives. The corresponding set of optimal solutions constitute the Pareto set. However, it has been pointed out that when using weights to reflect the relative preferences, effects of objectives magnitude differences should be avoided\cite{marler2004survey}. Hence, we apply $\delta_e$ and $\delta_d$ to normalize two objectives to a similar order of magnitude.
Constraint (\ref{con:a}) represents caching decisions at the SBSs. Constraint (\ref{con:b}) expresses the cache capacity of SBSs. 
Constraint (\ref{con:c}) enforces that the fraction of bandwidth not exceed the available bandwidth. 
Constraint (\ref{con:d}) is a QoS requirement, where $\gamma_i$ refers to data rate requirement of the $i^{th}$ user.
Constraint (\ref{con:e}) ensures that the backhaul link data rate is greater than a threshold.
We solve the proposed  optimization problem by using the MOSEK solver and CVX with the associated built-in functions on MATLAB2020.

\section{Numerical Investigations and Analysis}\label{results}

In this section, we investigate numerically the trade-off between energy consumption and delay of a typical cache-enabled HetNet as the one shown in Fig. \ref{fig:Network Topology}. The impacts of caching time, caching capacity, SBSs density and transmission power on the system performance are analysed as well. 

\subsubsection{Simulation Settings}
The simulation topology is composed of the MBS and SBSs in a 400m-by-400m square based on a random distribution, and one SBS serves its associated users within a 40m radius area. The number of UEs and popular content files are set to be 100 and 200, respectively. Without loss of generality the size of each file is uniformly set to be $10^7$ bits\cite{gabry2016energy}. The key simulation parameters assumed in this study are summarized in Table \ref{tab:parameter}. 

\begin{table}[ht!]
\caption{ Simulation Parameters\cite{gabry2016energy}\cite{xu2014coordinated}\cite{saha2019millimeter}}
\begin{center}
    \begin{tabular}{l|c}
    \hline
    \textbf{Parameter} &\textbf{value}\\
    \hline
    \hline
    $R_M$: radius of MBS coverage & 400m \\
    \hline
    $R_S$: radius of SBS coverage & 40m \\
     \hline
     $Q$: size of each file      & $10^6$ bits \\
     \hline
     $M$: number of files      &  200 \\
     \hline
     $N$: cache capacity per SBS &  50, 100, 150, 200 \\
     \hline
     $T$: caching time &  0.1, 1, 10, 100 hours\\
     \hline
    $f_c$: carrier frequency    & 28GHz \\
    \hline
    $W$: bandwidth of system.   & 200MHz \\
    \hline
    $w^{ca}$: caching power consumption per bit  & $6.25\cdot 10^{-12}$W/bit \\
    \hline
    $P^t_M$: transmit power of the MBS  & 46dBm\\
    \hline
    $P^t_S$:transmit power per SBS  & 20, 23, 26dBm\\
    \hline
    $G^a$: antenna gain   & 18dBi\\
    \hline
    $N_0$: noise power spectral efficiency  & -173dBm/Hz\\
    \hline
    $\delta^l$,$\delta^n$: path loss index for LOS and NLOS & $\delta^l$=2, $\delta^n$=3.3\\
    \hline
    \end{tabular}
    \label{tab:parameter}
    \end{center}
\end{table}

\subsubsection{Effect of caching time on the different schemes}

\begin{figure*}[htbp]
\centering
\begin{subfigure}{0.5\columnwidth}
\centering
\includegraphics[width=\columnwidth]{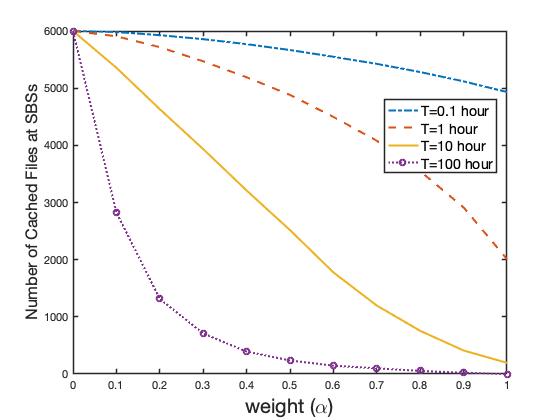}
\caption{caching decision}
\label{fig:timeca}
\end{subfigure}
\begin{subfigure}{0.5\columnwidth}
\centering
\includegraphics[width=\columnwidth]{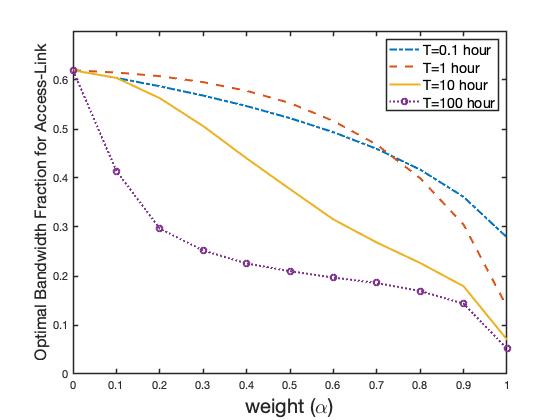}
\caption{bandwidth allocation}
\label{fig:timebw}
\end{subfigure}
\begin{subfigure}{0.5\columnwidth}
\centering
\includegraphics[width=\columnwidth]{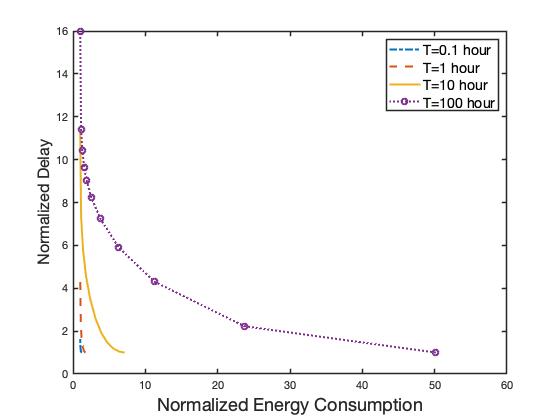}
\caption{Pareto front}
\label{fig:time pareto}
\end{subfigure}
\begin{subfigure}{0.5\columnwidth}
\centering
\includegraphics[width=\columnwidth]{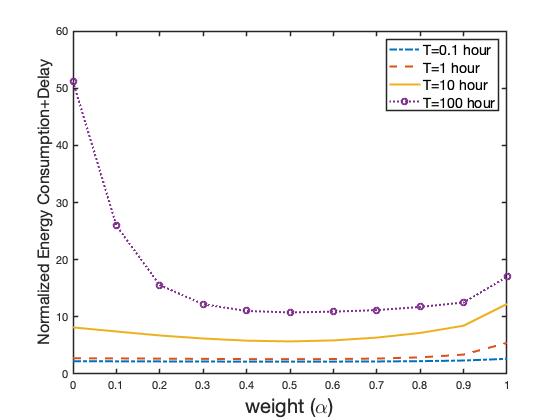}
\caption{sum of energy consumption and delay}
\label{fig:time objectivie}
\end{subfigure}
\caption{Simulation results in terms of caching time}
\label{fig:capacity}
\end{figure*}

This subsection evaluates the impact of the caching time on system performance. For the sake of comparison, we also present the baseline results of delay-only and energy-only optimization schemes in accordance to \cite{bacstuug2016delay}\cite{gabry2016energy}. Those schemes can be regraded as the case when the preference coefficient $\alpha$ equals to $0$ and $1$, respectively.

Undoubtedly, when compared to the actual transmission time, the caching stage takes place on a much larger scale; popular content might be updated in a daily or hourly basis\cite{tao2016content}. Hence, the caching time is related to the files' life-cycle, or the frequency with which cached files are renewed. If the requests for a content are concentrated in a short period, we say that the life-cycle is short. On the contrary, if the requests are spread out over a long time, the life-cycle is a comparatively long one. Therefore, in terms of different caching time, Fig. \ref{fig:time pareto} depicts the trade-off between energy consumption and delay under the same topology and normalization parameters. As shown in the Pareto curve, the energy consumption increases while the delay declines, indicating that the energy consumption and delay are indeed competitive. Moreover, it can be illustrated from Fig. \ref{fig:time objectivie}, the sum of energy consumption and delay at the extreme points are always higher than the rest of the points of same curve. For example, when $T=100 h$, comparing the lowest point where $\alpha=0.5$ with the points of the delay-only ($\alpha=0$) and the energy only ($\alpha=1$) schemes, the sum is decreased by $80\%$ and $50\%$, respectively, which indicates that the proposed scheme can achieve a better network performance.

Moreover, Fig. \ref{fig:timeca} and Fig. \ref{fig:timebw} show how edge caching decisions and bandwidth allocation vary with the weighting ($\alpha$). As the weighting of the energy consumption increases, fewer files are cached at SBSs to reduce the energy consumption and a larger fraction of bandwidth is assigned to the backhaul links to meet the corresponding increasing load. The delay increases as more users cannot fetch the files from SBSs directly. Aside from that, observe that even though the same capacity and bandwidth are assumed, the number of cached files and the fraction of bandwidth allocated to the access link generally decreases as caching time increases. Therefore, a dynamic resource allocation is necessary to attain the optimal results. Furthermore, in a practical scenario, short-lived popular content should be cached at the SBSs, while long-lived one should be anchored at the data center instead.

\subsubsection{Comparison for varying  caching capacity of SBS}

\begin{figure*}[htbp]
\centering
\begin{subfigure}{0.67\columnwidth}
\centering
\includegraphics[width=\columnwidth]{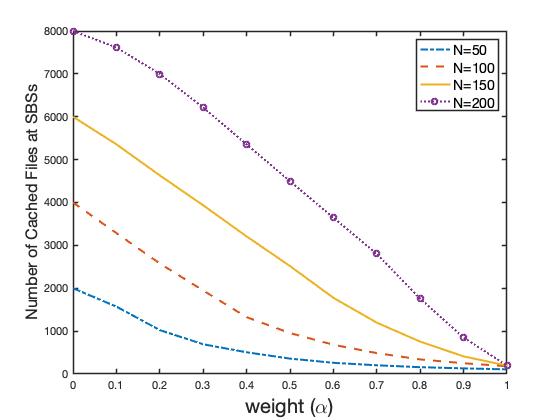}
\caption{caching decision}
\label{fig:capacity_ca}
\end{subfigure}
\begin{subfigure}{0.67\columnwidth}
\centering
\includegraphics[width=\columnwidth]{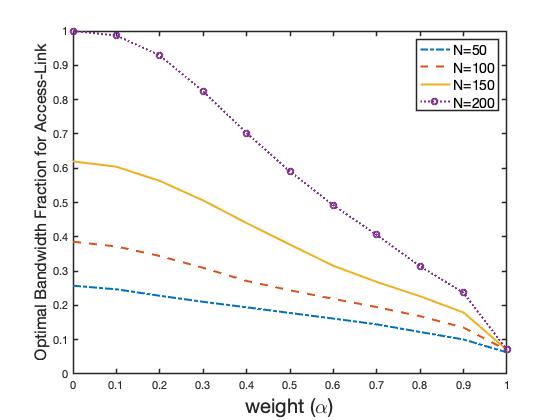}
\caption{bandwidth allocation}
\label{fig:capacity_bw}
\end{subfigure}
\begin{subfigure}{0.67\columnwidth}
\centering
\includegraphics[width=\columnwidth]{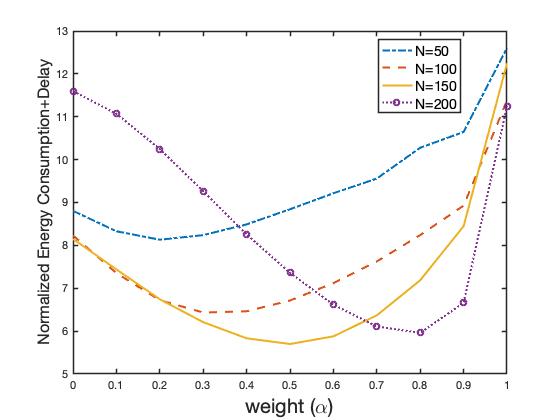}
\caption{sum of energy consumption and delay}
\label{fig:capacity_obj}
\end{subfigure}
\caption{Simulation results in terms of capacity}
\label{fig:capacity}
\end{figure*}

Hereafter we assess network performances in terms of cache capacity. For ease of analysis and without loss of generality, the cache capacity of SBSs is set to be the same. 

\begin{figure*}[htbp]
\centering
\begin{subfigure}{0.67\columnwidth}
\centering
\includegraphics[width=\columnwidth]{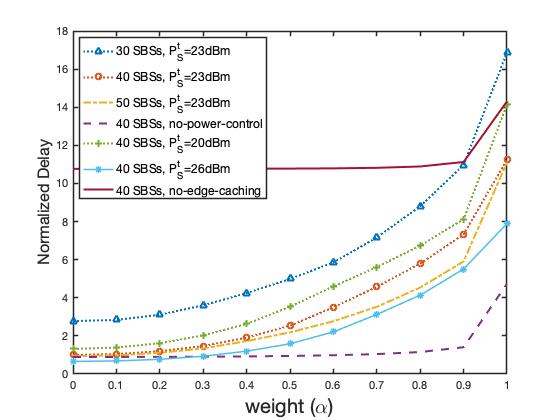}
\caption{aggregated delay}
\label{fig:density_power_de}
\end{subfigure}
\begin{subfigure}{0.67\columnwidth}
\centering
\includegraphics[width=\columnwidth]{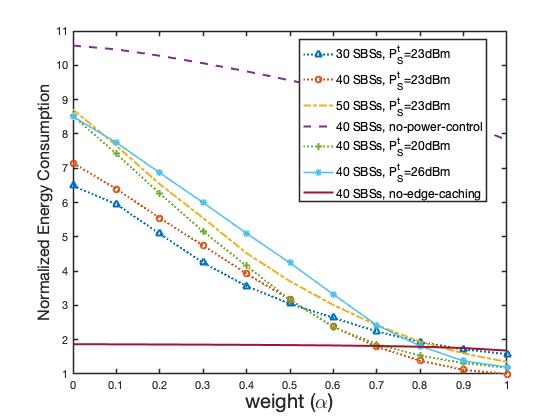}
\caption{energy consumption}
\label{fig:density_power_en}
\end{subfigure}
\begin{subfigure}{0.67\columnwidth}
\centering
\includegraphics[width=\columnwidth]{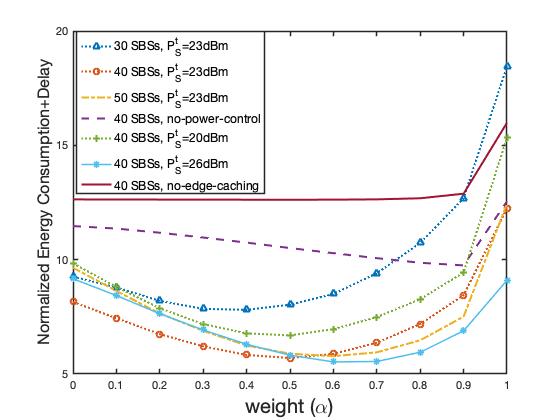}
\caption{sum of energy consumption and delay}
\label{fig:density_power_obj}
\end{subfigure}
\caption{Simulation results in terms of density and transmission power. }
\end{figure*}
As shown in Fig. \ref{fig:capacity_ca} and Fig. \ref{fig:capacity_bw}, the amount of files cached at SBSs and the bandwidth portion allocated to the access link decrease with the increased energy consumption preference. This is because reducing edge caching to lower caching energy consumption results in an increased traffic load being redirected to the backhaul link, which subsequently leads to more corresponding bandwidth being required. Furthermore, observe in Fig. \ref{fig:capacity_bw} that in the limiting case where the emphasis is only on energy consumption ($\alpha$ = 1), regardless of the available caching capacity, the number of cached files converges to a single value, inferring that any disparities due to capacity limits becoming less significant. 

The general trend depicted in Fig. \ref{fig:capacity_obj} is that the sum decreases with increased capacity, indicating that increasing capacity can improve network performance. However, note that when capacity is large enough, further increasing the available capacity might degrade the network. This is because when delay optimization is preferred, edge caching is prioritized which comes at a greater expense of caching energy consumption and is hence not necessarily a desirable operating point when considering both objectives jointly.

\subsubsection{Comparison in small cell density and available transmission power}

This subsection illustrates the results of different small cell densities and transmission power, furthermore, a scheme without edge caching in accordance to \cite{saha2019millimeter} is presented. Based on the above results, a caching life-cycle of 10 hours and a capacity of 150 files are assumed for the following analysis. 

As shown from Fig. \ref{fig:density_power_obj}, the overall performance is improved  with an increased density, owing to the fact that a denser set of SBSs can respond to the user's requests via a reduced distance, resulting in a lower delay. However, there is a saturation point where increasing further the density leads to very strong interference and higher energy consumption. As a result, an overall network degradation can be observed. Furthermore, when the network with different transmission power levels is considered, it is intuitive that the scheme without edge caching capability generates the worst performance due to the longest delays in fetching content to the end users. Observe also that when SBSs transmit at full power (no power control) the performance in terms of delay is very competitive, however, the energy consumption is three (3) times higher than the case when fractional power control is in place. Note, that if the maximum transmission power is halved to 20dBm will not enhance the performance either, because reducing transmission power can increase the latency by almost $40\%$. At this point, it is also worth observing that doubling the transmission power to 26dBm can outperform the case with the primal transmission power when $\alpha \textgreater 0.5$, but under-performs that for the rest of the cases. This is because, as 
presented in Fig. \ref{fig:density_power_de} and Fig.\ref{fig:density_power_en}, when $\alpha$ increases, optimizing energy consumption is prioritized, leading to the disparities in energy consumption narrowing and impact of latency difference is dominant on the overall value difference, which makes the case with a higher transmission power performing better. On the contrary, when $\alpha$ decreases, impact of energy consumption difference is dominant instead, which making the network with a lower transmission power performing better.

In summary, for the varying  power transmission investigation, the proposed framework allows to provide delay-energy consumption network operating points, which achieve a superior performance over single objective, albeit optimal, optimization schemes with a proper selection of the weighting parameter $\alpha$.
Taking $\alpha=0.5$ for example, when compared with the delay-only and the energy-only optimization schemes, the summation of energy consumption and delay can be reduced by  $30\%$ and $55\%$ respectively. 

\section{Conclusions}\label{conclusion}
In this paper,  the trade-off between energy consumption and delay in a cache-enabled integrated access and backaul (IAB) mmWave HetNet via a weighted-sum multi-objective optimization problem, where popular content cache placement and mmWave bandwidth partitioning is considered jointly. More specifically, we provide a comprehensive view on the impact of a number of key parameters such as the caching time, the small cell capacity and density and the transmission power on the system performance. A wide set of numerical investigations reveal that when caching is enabled by simply upgrading network capacity, density of small cells or the maximum transmission power may not always improve network performance. Another important observation is that the actual life-cycle management of popular content need to be considered when caching at the edge is selected. We have shown that a desirable Pareto efficient point corresponding to  the minimum of the aggregation of the two objectives can be obtained with an appropriate weighting coefficient, resulting in significant performance improvement compared to the case where delay to access popular content or energy consumption are optimized independently.

\bibliographystyle{IEEEtran}
\bibliography{IEEEabrv,EDPareto}

\begin{thebibliography}{10}
\providecommand{\url}[1]{#1}
\csname url@samestyle\endcsname
\providecommand{\newblock}{\relax}
\providecommand{\bibinfo}[2]{#2}
\providecommand{\BIBentrySTDinterwordspacing}{\spaceskip=0pt\relax}
\providecommand{\BIBentryALTinterwordstretchfactor}{4}
\providecommand{\BIBentryALTinterwordspacing}{\spaceskip=\fontdimen2\font plus
\BIBentryALTinterwordstretchfactor\fontdimen3\font minus
  \fontdimen4\font\relax}
\providecommand{\BIBforeignlanguage}[2]{{%
\expandafter\ifx\csname l@#1\endcsname\relax
\typeout{** WARNING: IEEEtran.bst: No hyphenation pattern has been}%
\typeout{** loaded for the language `#1'. Using the pattern for}%
\typeout{** the default language instead.}%
\else
\language=\csname l@#1\endcsname
\fi
#2}}
\providecommand{\BIBdecl}{\relax}
\BIBdecl

\bibitem{lv2017interference}
W.~Lv, Z.~Zhang, C.~Jiao, and C.~Zhong, ``Interference coordination in
  full-duplex hetnet with large-scale antenna arrays,'' in \emph{2017 IEEE
  ICC}.\hskip 1em plus 0.5em minus 0.4em\relax IEEE, 2017, pp. 1--6.

\bibitem{bacstuug2016delay}
E.~Ba{\c{s}}tu{\u{g}}, M.~Kountouris, M.~Bennis, and M.~Debbah, ``On the delay
  of geographical caching methods in two-tiered heterogeneous networks,'' in
  \emph{2016 IEEE 17th SPAWC}.\hskip 1em plus 0.5em minus 0.4em\relax IEEE,
  2016, pp. 1--5.

\bibitem{chandrasekhar2008femtocell}
V.~Chandrasekhar, J.~G. Andrews, and A.~Gatherer, ``Femtocell networks: a
  survey,'' \emph{IEEE Comm. mag.}, vol.~46, no.~9, pp. 59--67, 2008.

\bibitem{khan2019energy}
B.~S. Khan, S.~Jangsher, H.~K. Qureshi, and S.~Mumtaz, ``Energy efficient
  caching in cooperative small cell network,'' in \emph{2019 16th IEEE
  CCNC}.\hskip 1em plus 0.5em minus 0.4em\relax IEEE, 2019, pp. 1--6.

\bibitem{wu2020energy}
H.~Wu, H.~Lu, F.~Wu, and C.~W. Chen, ``Energy and delay optimization for
  cache-enabled dense small cell networks,'' \emph{IEEE Transactions on
  Vehicular Technology}, vol.~69, no.~7, pp. 7663--7678, 2020.

\bibitem{song2021energy}
J.~Song, Q.~Song, Y.~Wang, and P.~Lin, ``Energy-delay tradeoff in adaptive
  cooperative caching for energy-harvesting ultradense networks,'' \emph{IEEE
  Transactions on Computational Social Systems}, 2021.

\bibitem{gabry2016energy}
F.~Gabry, V.~Bioglio, and I.~Land, ``On energy-efficient edge caching in
  heterogeneous networks,'' \emph{IEEE JSAC}, vol.~34, no.~12, pp. 3288--3298,
  2016.

\bibitem{wang2019energy}
L.~Wang, H.~Zhang, J.~Qiao, X.~Zhou, and D.~Yuan, ``Energy-delay aware user
  association in mmwave backhaul networks using matching theory,'' in \emph{ICC
  2019-2019 IEEE ICC}.\hskip 1em plus 0.5em minus 0.4em\relax IEEE, 2019, pp.
  1--6.

\bibitem{zhang2021survey}
Y.~Zhang, M.~A. Kishk, and M.-S. Alouini, ``A survey on integrated access and
  backhaul networks,'' \emph{arXiv preprint arXiv:2101.01286}, 2021.

\bibitem{saha2019millimeter}
C.~Saha and H.~S. Dhillon, ``Millimeter wave integrated access and backhaul in
  5g: Performance analysis and design insights,'' \emph{IEEE JSAC}, vol.~37,
  no.~12, pp. 2669--2684, 2019.

\bibitem{bai2014coverage}
T.~Bai and R.~W. Heath, ``Coverage and rate analysis for millimeter-wave
  cellular networks,'' \emph{IEEE Transactions on Wireless Communications},
  vol.~14, no.~2, pp. 1100--1114, 2014.

\bibitem{xu2014coordinated}
Y.~Xu, Y.~Li, Z.~Wang, T.~Lin, G.~Zhang, and S.~Ci, ``Coordinated caching model
  for minimizing energy consumption in radio access network,'' in \emph{2014
  IEEE ICC}.\hskip 1em plus 0.5em minus 0.4em\relax IEEE, 2014, pp. 2406--2411.

\bibitem{jindal2008fractional}
N.~Jindal, S.~Weber, and J.~G. Andrews, ``Fractional power control for
  decentralized wireless networks,'' \emph{IEEE Transactions on Wireless
  Communications}, vol.~7, no.~12, pp. 5482--5492, 2008.

\bibitem{marler2004survey}
R.~T. Marler and J.~S. Arora, ``Survey of multi-objective optimization methods
  for engineering,'' \emph{Structural and multidisciplinary optimization},
  vol.~26, no.~6, pp. 369--395, 2004.

\bibitem{tao2016content}
M.~Tao, E.~Chen, H.~Zhou, and W.~Yu, ``Content-centric sparse multicast
  beamforming for cache-enabled cloud ran,'' \emph{IEEE Transactions on
  Wireless Communications}, vol.~15, no.~9, pp. 6118--6131, 2016.

\end{thebibliography}
\end{document}